\documentclass[11pt]{amsart}
\usepackage[T1]{fontenc}
\usepackage[utf8]{inputenc}
\usepackage{lmodern}
\usepackage{amsmath,amssymb,amsthm,mathtools,mathrsfs}
\usepackage[a4paper,margin=2.8cm]{geometry}
\usepackage[hypertexnames=false,hidelinks]{hyperref}
\numberwithin{equation}{section}
\theoremstyle{plain}
\newtheorem{proposition}{Proposition}[section]
\allowdisplaybreaks[1]
\title[Old structures, a vanishing formula]
{Old Structures, a Vanishing Formula}
\author{K. Castillo}
\address{CMUC, Department of Mathematics, University of Coimbra,
3000-143 Coimbra, Portugal}
\email{math@keniercastillo.com}
\author{G. Gordillo-N\'u\~nez}
\address{CMUC, Department of Mathematics, University of Coimbra,
3000-143 Coimbra, Portugal}
\email{up202310693@up.pt}
\date{11 October 2026}
\hypersetup{
  pdftitle={Old Structures, a Vanishing Formula},
  pdfauthor={K. Castillo and G. Gordillo-N\'u\~nez},
  pdfsubject={Comparison of fourth-order Hahn constructions with the 1998 elimination},
  pdfkeywords={Hahn operator, orthogonal polynomials, fourth-order equation, normalization}
}
\begin{document}
\begin{abstract}
We identify a normalization defect in a determinantal construction
of fourth-order Hahn difference equations and compare its scalar
equation with Foupouagnigni's 1998 elimination and its published
2001 counterpart.
On an explicit positive orthogonal family, all five prescribed
coefficients vanish: the formula gives \(0=0\) and imposes no
condition. Along an explicit regular perturbation, division by
\(\varepsilon^2\) before specialization recovers a nonzero equation
on that family, the same scalar equation as the correctly executed
old shift-and-eliminate procedure. On the common nondegenerate domain specified below, both
constructions give the same scalar operator up to a nonzero
rational factor. The first-order
degree bounds follow directly from valid identities of 1998;
preserving cancellations gives a bound that can be arbitrarily
sharper on explicit admissible inputs, together with exact top
coefficients and cancellation criteria. Explicit coefficient formulas
do not by themselves yield a new scalar equation on this domain.
Our comparison does not exclude corrective or expository contributions;
additional computational gain or exceptional coverage must be demonstrated.
\end{abstract}
\maketitle

\section{The equation on the common nondegenerate domain}

Theorem~4 of~\cite{MMMRP2026} prescribes fourth-order Hahn
coefficients. Proposition~\ref{prop:zero-normalization} exhibits a positive family on which
they all vanish, giving \(0=0\), and a normalization recovering
the scalar equation yielded by the 1998 elimination. Proposition~\ref{prop:known} identifies the
two scalar operators, up to a nonzero rational factor, on their
common nondegenerate domain; the old procedure is
\cite[Theorem~6.4, p.~87, Eq.~(6.53)]{F1998}.
Over \(\mathbb C(y)\), a nonzero factor preserves an equation;
a zero factor erases it. No classification of every exceptional
input is asserted. Thesis citations use printed page numbers.

Fix \(0<q<1\) and \(\omega>0\). Put
\(c=\omega/(1-q)\), \(y=x-c\), \(s=(q-1)y\), and define
\[
 Tf(y)=f(qy),\quad \Delta=(T-I)/s.
\]
Then \(qx+\omega-c=qy\): the Hahn operator becomes Jackson's
operator by exact translation.
The class equivalence is Theorem~4.2(i), pp.~67--68, of~\cite{F1998};
the polynomial transport is its part~(ii), Eq.~(4.21), p.~68,
and the operator conjugacy is Eq.~(4.23) on the same page.
Let \(\mathcal U:\mathbb C[y]\to\mathbb C\) be linear, with \(\mathcal U[1]=1\)
and \(\det(\mathcal U[y^{j+k}])_{j,k=0}^{r}\ne0\) for every \(r\geq0\).
Its formal Stieltjes series is
\(S(y)=\sum_{j\geq0}\mathcal U[y^j]y^{-j-1}\in\mathbb C((y^{-1}))\).
Use its monic orthogonal and first associated sequences \(P_n\),
\(R_n=P_n^{(1)}\), with \(P_{-1}=R_{-1}=0\), \(P_0=R_0=1\),
recurrence coefficients \(\beta_n\) for \(n\geq0\), \(\gamma_n\ne0\) for \(n\geq1\),
and \(\gamma_0=1\). The polynomial data \(a,B,C_1,C_2,D\in\mathbb C[y]\)
in the transported Riccati identity refer to this normalized series:
\begin{equation}
 a\Delta S=B S(TS)+C_1S+C_2(TS)+D,\quad a,B\not\equiv0,
 \label{eq:riccati}
\end{equation}
Here \(a(y)=A(c+qy)\), with the other coefficients transported likewise.
Products and orders at infinity are formal. If \(S_{\mathrm{raw}}=\sigma S\),
\(\sigma\ne0\), normalization gives
\((B,D)=(\sigma B_{\mathrm{raw}},D_{\mathrm{raw}}/\sigma)\),
leaving the other Riccati coefficients unchanged.

For a fixed \(n\geq1\), set
\begin{equation}
 \begin{aligned}
 \Phi_n&=\begin{pmatrix}P_{n+1}&R_n\\[7pt]P_n&R_{n-1}\end{pmatrix},
 \quad p=\Theta_{n-1}/\gamma_n,\quad
 \ell=L_{n-1}+(y-\beta_n)p,\\[7pt]
 M_n&=\begin{pmatrix}L_n&\Theta_n\\[7pt]-p&\ell\end{pmatrix},
 \quad C=\begin{pmatrix}C_1&-D\\[7pt]B&-C_2\end{pmatrix}.
 \end{aligned}
 \label{eq:matrices}
\end{equation}
Equations~(42)--(43) of~\cite{MMMRP2026} are
\begin{equation}
 (aI+sM_n)\Phi_n=(T\Phi_n)(aI+sC).
 \label{eq:new-system}
\end{equation}

\begin{proposition}
\label{prop:known}
Compare degree \(m=n+1\) in Theorem~4 of~\cite{MMMRP2026}
and association order zero in the elimination of
\cite[Theorem~6.3, p.~84; Theorem~6.4, p.~87]{F1998}.
Suppose the polynomial
coordinate is cyclic of rank four in their structure system,
and both scalar constructions have nonzero fourth-order leading
coefficients. Their monic normalizations in powers of \(T\)
are identical over \(\mathbb C(y)\).
\end{proposition}

\begin{proof}
The determinant of \(\Phi_n\) is a nonzero constant, by the
three-term recurrence. Taking determinants in
\eqref{eq:new-system} therefore gives the coefficient identity
\[
 \det(aI+sM_n)=\det(aI+sC)=:\delta.
\]
Regularity excludes rational \(S\), whose Hankel rank would be
finite. In the rearranged identity
\[
 TS=\frac{(a+sC_1)S+sD}{-sBS+a-sC_2}
\]
the denominator cannot vanish, since \(B\not\equiv0\).
If \(\delta=0\), the fractional linear map is constant over
\(\mathbb C(y)\), making \(TS\), hence \(S\), rational.
Thus \(\delta\not\equiv0\).

For a two-by-two matrix \(X\),
\(\operatorname{adj}(aI+sX)=aI+s\operatorname{adj}X\).
Adjugation of \eqref{eq:new-system} thus gives
\begin{equation}
 a\Delta\Phi_n=U_n(T\Phi_n)+\Phi_n V,\quad
 U_n=-\operatorname{adj}M_n,\quad V=\operatorname{adj}C.
 \label{eq:old-system}
\end{equation}
The thesis's sign \(-S\)
\cite[Definition~4.1, p.~61, Eq.~(4.1)]{F1998} and \(\gamma_0=1\)
give the same \(V\) in its Riccati equation (6.1), p.~77.
Its Theorem~6.2, p.~82, identifies the left matrix with \(U_n\),
since \(T\Phi_n\) is invertible. Hence \eqref{eq:old-system}
is its (6.29)--(6.32), with dictionary
\begin{equation}
 \begin{aligned}
 (G_0,E_0,F_0,H_0)&=(-B,C_1,C_2,-D),\\[7pt]
 (E_{n+1},F_{n+1},H_{n+1})&=(L_n,-\ell,-\Theta_n),\quad
 H_n/\gamma_n=-p.
 \end{aligned}
 \label{eq:dictionary}
\end{equation}
The same four relations are published in
\cite[Lemmas~2--3, p.~453, Eqs.~(25)--(28)]{FHR2001}.
The identity \(F_{n+1}=-\ell\) also agrees with its recurrence
(6.5), p.~78. Since \(\delta\ne0\), adjugation is reversible:
the two systems define the same linear shift map for
\emph{every formal matrix state}, not only for \(\Phi_n\).

At degree \(m=n+1\), the old elimination uses the lower-left
coordinate at index \(m\), the new one the upper-left at
\(m-1\). The recurrence gauge is
\[
 \Phi_m=W_m\Phi_{m-1},\quad
 W_m=\begin{pmatrix}y-\beta_m&-\gamma_m\\[7pt]1&0\end{pmatrix}.
\]
Since \(\det W_m=\gamma_m\ne0\), compatibility (45)
in~\cite{MMMRP2026},
\(a\Delta W_m=M_mW_m-(TW_m)M_{m-1}\),
and the product rule show that \(X\mapsto W_mX\) gauges every
formal state from index \(m-1\) to \(m\). Its lower-left
coordinate is the old upper-left one, so both eliminations give
universal scalar relations for the same observable. For 1998 we
execute (6.51)--(6.53), p.~87, rather than copying the expanded
list.\footnote{The list below (6.52), p.~87, prints
\(b_2(Tb_2)(Ta_1)\) in the \(T^3\) coefficient instead of
\(b_2(Tb_2)(T^2a_1)\). The publication gives the shift instructions
\cite[Eqs.~(43)--(44), p.~455]{FHR2001} without that list.
The pairs with coefficients (6.40), p.~84, follow from (6.43)--(6.46),
p.~85, by multiplication and substitution, without dividing by
\(K_3\) or \(K_7\); hence they remain valid for every formal state
when an intermediate divisor vanishes.}

Vectorize the common system as \(Tz=H(y)z\), with polynomial
coordinate \(f=(1,0,0,0)z\). Set
\begin{equation}
 h_0=(1,0,0,0),\quad h_{j+1}=(Th_j)H,\quad
 \mathcal C=\begin{pmatrix}h_0\\[7pt]h_1\\[7pt]h_2\\[7pt]h_3\end{pmatrix}.
 \label{eq:cyclic}
\end{equation}
Since \(T^jf=h_jz\) and \(\det\mathcal C\not\equiv0\), the
first four observations determine every state, giving the unique
monic universal relation
\[
 L_*=T^4+\sum_{j=0}^3b_jT^j,\quad
 (b_0,b_1,b_2,b_3)=-h_4\mathcal C^{-1}.
\]
Two distinct such relations would differ by one of order at most three,
contradicting independence of \(h_0,\ldots,h_3\). For their
nonzero leading coefficients \(a_4,b_4\) in powers of \(T\),
\begin{equation}
 \frac{1}{a_4}\mathcal O_{26}=\frac{1}{b_4}\mathcal O_{98}=L_*.
 \label{eq:equivalence}
\end{equation}
Coefficients multiply from the left. This identifies entire
operators, not all equations annihilating one polynomial.
Triangular conversion between \(T\) and \(\Delta\), with
nonzero leading factors, preserves the equivalence.
\end{proof}

Lower cyclic rank is outside Proposition~\ref{prop:known}:
fourth-order multiples of a lower-order relation need not be
proportional.

\section{First-order bounds from the old identities}
\label{sec:bounds}

The degree estimates in Theorem~1 of~\cite{MMMRP2026} are
elementary consequences of the 1998 structure and remainder
identities. Keeping cancellations before counting degrees gives
the following refinement. Here \(B\) may vanish, as in that
theorem, and \(\deg0=-\infty\).

\begin{proposition}
\label{prop:degrees}
For a normalized regular input satisfying the Riccati identity
\eqref{eq:riccati}, with \(a\not\equiv0\) and \(B\) allowed to
vanish, put
\[
 u=a-sC_2,\quad v=C_1+C_2,\quad
 t=\max\{\deg u-2,\deg B-2,\deg v-1\}.
\]
For every \(n\geq0\), the polynomial structure coefficients are
uniquely determined by the fixed input and satisfy
\begin{equation}
 \begin{aligned}
 \deg\Theta_n&\leq t,\\[7pt]
 \deg(L_n-C_1)&\leq
 \max\{\deg(a+sC_1)-1,\deg B-1,\deg\Theta_n-1\},\\[7pt]
 \deg(L_n+C_2)&\leq
 \max\{\deg u-1,\deg D+1,\deg\Theta_n-1\}.
 \end{aligned}
 \label{eq:degree-bounds}
\end{equation}
If \(t\geq0\), write \([y^j]f\) for the coefficient of \(y^j\)
in \(f\), put \([j]_q=(q^j-1)/(q-1)\),
\(\alpha=[y^{t+2}]u\), \(b=[y^{t+2}]B\),
\(e=[y^{t+1}]v\). Then
\begin{equation}
 \begin{aligned}
 [y^t]\Theta_n
 &=-\gamma_{n+1}q^{-n-2}
 \big\{\alpha[2n+3]_q+b(1+q^{2n+2})+eq^{2n+3}\big\},\\[7pt]
 [y^{t+1}](L_n+C_2)
 &=\alpha[n+1]_q+q^n(b+qe).
 \end{aligned}
 \label{eq:top-coefficients}
\end{equation}
The corresponding degrees are \(t\) and \(t+1\) exactly
when these coefficients are nonzero; otherwise they are at most
\(t-1\) and \(t\), respectively. Each top coefficient can vanish
at at most one index, unless it vanishes for every \(n\).
Uniform cancellation of the top coefficient is equivalent to
\(b=\alpha/(q-1)\), \(e=-(1+q^{-1})b\) for \(\Theta_n\),
and to \(\alpha=0\), \(b+qe=0\) for \(L_n+C_2\).
\end{proposition}

\begin{proof}
Write \(P=P_{n+1}\), \(R=R_n\), and
\(\kappa=RP_n-PR_{n-1}=\prod_{j=1}^{n}\gamma_j\ne0\),
with \(\kappa=1\) at \(n=0\).
The old structure identities
\cite[Theorem~6.2, p.~82, Eqs.~(6.29)--(6.32)]{F1998}
give the coefficient dictionary \eqref{eq:dictionary} by
adjugation, which requires no division by the transfer
determinant for the actual polynomial matrix. The constant
Casoratian \(\kappa\) makes the pair \((L_n,\Theta_n)\) unique.
Initial recurrences (6.4), (6.6), p.~78, give the same
identification at \(n=0\).
Set \(\eta=S-R/P\). The tail and remainder identities
\cite[Lemmas~6.1--6.2, p.~81, Eqs.~(6.21)--(6.22)]{F1998}, together with the
Casoratian (2.24), p.~23, give
\[
 \eta=\kappa\gamma_{n+1}y^{-2n-3}+O(y^{-2n-4}),\quad
 \Delta\eta=O(y^{-2n-4}).
\]
The old coefficient residual (6.27), p.~82, becomes
\begin{equation}
 \kappa\Theta_n=P(TP)
 \big[u\Delta\eta-v\eta
 -B\{S(T\eta)+\eta(TS)-\eta(T\eta)\}\big].
 \label{eq:degree-residual}
\end{equation}
This uses only \(T\eta=\eta+s\Delta\eta\) to collect the
linear terms. Since \(S=O(y^{-1})\) and
\(\deg P(TP)=2n+2\), its three groups have orders at most
\(\deg u-2\), \(\deg v-1\), and \(\deg B-2\).
The left side is polynomial, proving the first bound.
The other two follow by comparing degrees in the old
structures, rewritten as
\[
 \begin{aligned}
 (a+sC_1)\Delta P&=(L_n-C_1)P+\Theta_nP_n-B(TR),\\[7pt]
 u\Delta R&=(L_n+C_2)R+\Theta_nR_{n-1}+D(TP).
 \end{aligned}
\]
At \(n=0\), \(R=1\) and \(R_{-1}=0\), so the same bounds
hold.
The reduced input \(u\Delta S=B S(TS)+vS+D\) gives
\(\deg D\leq t\), hence \(\deg(L_n+C_2)\leq t+1\).
Leading terms in the residual and the first structure give
\eqref{eq:top-coefficients}; after removing nonzero prefactors,
these are affine in \(q^{2n+3}\) and \(q^n\), respectively.
Since \(0<q<1\), the cancellation criteria follow.
\end{proof}

Since \(\deg u\leq\max\{\deg a,\deg C_2+1\}\) and
\(\deg v\leq\max\{\deg C_1,\deg C_2\}\), \(t\) is no larger
than (33) of~\cite{MMMRP2026}. Adding \(C_1\) or subtracting
\(C_2\) in the shifted estimates recovers the minimum (34).
Both printed bounds are correct; they add no new first-order
construction.
We derive them from valid 1998 identities, without using its
different Proposition~6.2, p.~79, or the corresponding published
bounds \cite[Proposition~1(3), p.~450, Eq.~(13)]{FHR2001};
we do not claim their precise form was printed there.
Elementary derivation from valid old identities does not exclude
a corrective contribution relative to older degree estimates.
For classical Pearson data with \(B=C_2=0\),
\eqref{eq:top-coefficients} specializes to the constants already
printed in \cite[p.~458, Eq.~(52)]{FHR2001},
with association index \(n+1\), under dictionary~\eqref{eq:dictionary}.

For \(u\ne0\), applying Theorem~1 to the reduced input with
\(C_2=0\) already gives \(t\); our residual also covers
\(u=0\). The refinement removes redundancy.

The improvement can be arbitrarily large under the stated
hypotheses. On the positive family \eqref{eq:family}, with
data \eqref{eq:example-input}, replace
\[
 (a,C_1,C_2)\longmapsto
 (a+sy^k,C_1-y^k,C_2+y^k),\quad k\geq2.
\]
This preserves the Riccati identity, \(u\), \(v\) and
\(\Theta_n=-q^n\), while replacing \(L_n\) by \(1-y^k\).
There is still no common input factor, since \(B=D=-1\).
Estimate (33) gives only \(\deg\Theta_n\leq k-1\);
\eqref{eq:degree-bounds} gives the exact degree zero.
The shifted bounds also force \(L_n-C_1=0\) and give
\(\deg(L_n+C_2)\leq1\), attained by \(qy\).
The inflated estimate measures redundant presentation terms,
not additional polynomial structure.
The unmodified input \eqref{eq:example-input} already attains
both original bounds with \(\deg\Theta_n=\deg L_n=0\);
strict improvement on every admissible input is impossible.

\section{The zero formula and its normalized continuation}

Consider the monic family
\begin{equation}
 Q_0=1,\quad Q_1=y,\quad
 Q_{k+1}=yQ_k-q^{2k-1}Q_{k-1},\quad k\geq1.
 \label{eq:family}
\end{equation}
The bounded selfadjoint Jacobi matrix \(J\), with zero diagonal and
off-diagonal entries \(a_k=q^{k-1/2}>0\), has spectral measure
\(\mu\) at \(e_0\). The recurrence gives
\[
 Q_k(J)e_0=(a_1\cdots a_k)e_k,\quad
 \int Q_kQ_j d\mu=q^{k^2}\delta_{k,j}.
\]
These positive norms make all moment matrices nonsingular. With
\(F(y)=\int d\mu(t)/(y-t)\), deletion of the first Jacobi row
and column leaves \(qJ\). Its Schur complement yields
\[
 F(qy)=\frac1{qy-F(y)}.
\]
Consequently \eqref{eq:riccati} holds with
\begin{equation}
 a=-s,\quad B=D=-1,\quad C_1=1,\quad C_2=qy-1.
 \label{eq:example-input}
\end{equation}
Here \(a,B\ne0\), and there is no common input factor.
Translation gives an input with \(\omega>0\) in the original
Hahn variable. The associated polynomials are
\(R_k=q^kQ_k(y/q)\), and
\[
 TR_k=q^kQ_k,\quad
 TQ_{n+1}=q^{n+1}yQ_n-q^nR_{n-1}.
\]
Substitution in the structure relations gives
\(L_n=1\), \(\Theta_n=-q^n\).

Fix \(r=q^n\), \(d=q-1\). The two auxiliary coefficient
columns in the determinant of~\cite[(71)--(85)]{MMMRP2026}
have rows \(g_j\) satisfying
\begin{equation}
 \begin{aligned}
 g_1&=(-q^3r^2dy^2,\ q^3dy^3),\quad
 g_{j+1}=h\Delta g_j+(Tg_j)K,\\[7pt]
 \rho&=q^8rd^3y^5,\quad h=-dy\rho,\quad
 K=\begin{pmatrix}
 q^7d^3(qr-1)y^5&0\\[7pt]
 -q^6r^2d^3y^4&q^8rd^3y^5
 \end{pmatrix}.
 \end{aligned}
 \label{eq:rows}
\end{equation}
Multiplication gives \(g_2=\rho g_1\); the product rule and
\(h\Delta\rho+(T\rho)\rho=\rho^2\) then give
\(g_3=\rho^2g_1\). All three two-by-two minors vanish, and
each coefficient in Appendix~A is a linear combination of those
minors. Thus the printed operator is
\begin{equation}
 \mathcal O_0=\sum_{k=0}^4c_k\Delta^k=0,\quad
 c_0=c_1=c_2=c_3=c_4=0.
 \label{eq:zero}
\end{equation}
For any polynomial \(f\), the augmented determinant rows also
satisfy, by substitution in (79), (81), (83), (85),
\begin{equation}
 (g_j,F_j[f])=\rho^{j-1}(g_1,q^3rdy^2f),\quad j=1,2,3.
 \label{eq:augmented}
\end{equation}
The relations repeat the same information: eliminating the
auxiliaries also eliminates the condition on \(f\).
Thus \(0=0\) is a true identity imposing no restriction:
every polynomial satisfies it, regardless of degree or orthogonality.
Calling this a fourth-order equation would build
a castle in the air: the required nonzero operator is missing.

The order of the operations matters. In \(\varepsilon^2z=0\),
setting \(\varepsilon=0\) first gives \(0=0\); dividing by
\(\varepsilon^2\) for \(\varepsilon\ne0\) first gives \(z=0\),
which remains a condition at the limit. Likewise, our repair uses
the perturbed family before specialization; no finite multiple of
the already zero operator can recover its lost condition.
Common-factor cancellation is already discussed in
\cite[Remark~4, p.~456; Section~4.3, p.~462]{FHR2001}.
For real \(\varepsilon\), let
\(J_\varepsilon=J+\varepsilon e_0e_0^*\). It still has a positive
regular spectral measure. Its Stieltjes function and polynomials are
\[
 F_\varepsilon=\frac{F}{1-\varepsilon F},\quad
 Q_{\varepsilon,k}=Q_k-\varepsilon R_{k-1},
\]
where \(R_{-1}=0\), and \(R_k\) is unchanged.
The polynomial transform is \cite[p.~464, Eq.~(72)]{FHR2001}.
The new Riccati data are
\begin{equation}
 a=-s,\quad B=\varepsilon qy-\varepsilon^2-1,\quad
 C_1=1-\varepsilon,\quad C_2=qy-\varepsilon-1,\quad D=-1.
 \label{eq:perturbation}
\end{equation}
Their structure coefficients remain
\(L_n=1\), \(\Theta_n=-r\), \(\ell=1-y/r\).
Denote the prescribed operator by \(\mathcal O_\varepsilon\).
Equation~\eqref{eq:augmented} shows that every augmented
coefficient matrix has rank at most one at \(\varepsilon=0\).
Subtracting \(\rho\) and \(\rho^2\) times its first row from
the next two leaves two rows divisible by \(\varepsilon\).
Hence all five coefficients of \(\mathcal O_\varepsilon\) are
divisible by \(\varepsilon^2\).

The leading coefficient has exactly that order of vanishing.
In the notation of Appendix~A, \(c_4=I_3f_{3,4}\).
Direct expansion of its defining rows gives
\begin{equation}
 \begin{aligned}
 I_3&=\varepsilon q^{12}r^2d^5y^9(q^4r^2-1)+O(\varepsilon^2),\\[7pt]
 f_{3,4}&=\varepsilon q^{30}r^3d^{11}y^{15}+O(\varepsilon^2),\\[7pt]
 \left.\frac{c_4}{\varepsilon^2}\right|_{\varepsilon=0}
 &=q^{42}r^5d^{16}y^{24}(q^4r^2-1)\not\equiv0.
 \end{aligned}
 \label{eq:leading-limit}
\end{equation}
The shift system is regular in \(\varepsilon\), since
\(\det(aI+sC_\varepsilon)=s^2\). In the coordinate order
\(Q_{n+1},Q_n,R_n,R_{n-1}\), its transfer matrix at zero is
\[
 H_0=\begin{pmatrix}
 0&qry&0&-r\\[7pt]
 -qy/r&qy^2/r&1/r&-y/r\\[7pt]
 0&r&0&0\\[7pt]
 -1/r&y/r&0&0
 \end{pmatrix},\quad
 \det\mathcal C_0=qy^2(q^4r^2-1)\not\equiv0.
\]
Thus cyclic rank remains four; the example lies outside
Proposition~\ref{prop:known} because the printed operator
\(\mathcal O_0\) vanishes.
The system's unique monic fourth-order relation from
\eqref{eq:cyclic} is
\begin{equation}
 \begin{aligned}
 L={}&T^4-\frac{q^6}{r}y^2T^3
 +\left\{\left(q^7+\frac{q^3}{r^2}\right)y^2-1-q^2\right\}T^2\\[7pt]
 &-\frac{q^4}{r}y^2T+q^2I.
 \end{aligned}
 \label{eq:old-L}
\end{equation}
The 1998 elimination
\cite[Eq.~(6.40), p.~84; Eqs.~(6.51)--(6.53), p.~87]{F1998}
and the published determinant
\cite[Theorem~1, p.~456, Eq.~(45)]{FHR2001},
at \(m=n+1\), give the nonzero operator
\(q^{3-8m}(1-q^{2m+2})y^2L\).
The common paired-system tuple, (6.42), p.~84, in the thesis
and (33), p.~454, in the publication, is
\((0,-qy,q^{-m},1,0,-1,-q^{-m}y,-q^m)\);
after two shifts, the auxiliary rows of orders four and two are
proportional by \(q^{2-2m}\). Their subtraction gives
\eqref{eq:old-L}.

The regular shift system makes the divided limit a universal relation.
Cyclic uniqueness and \eqref{eq:leading-limit}, using the
\(\Delta^4\) coefficient \(q^6d^4y^4\) of \(L\), prove
\begin{equation}
 \left.\frac{\mathcal O_\varepsilon}{\varepsilon^2}
 \right|_{\varepsilon=0}
 =q^{36}r^5d^{12}y^{20}(q^4r^2-1)L\ne0.
 \label{eq:normalized-limit}
\end{equation}
The preceding computations prove the following statement.

\begin{proposition}
\label{prop:zero-normalization}
Fix \(n\geq1\), and put \(r=q^n\), \(d=q-1\).
The family \eqref{eq:family} is positive and regular, and the
data \eqref{eq:example-input} satisfy the hypotheses of
Theorem~4 of~\cite{MMMRP2026}, with no common input factor.
Nevertheless, all five prescribed coefficients vanish:
\(\mathcal O_0=0\). The polynomial coordinate of the structure
system is cyclic of rank four, and its unique monic universal
fourth-order relation is given by \(L\) in \eqref{eq:old-L}.

For every real \(\varepsilon\), the perturbation
\(J_\varepsilon=J+\varepsilon e_0e_0^*\) has a positive regular
spectral measure and Riccati data \eqref{eq:perturbation}.
Every coefficient of its prescribed operator
\(\mathcal O_\varepsilon\) is divisible by \(\varepsilon^2\),
and the leading coefficient has exactly that order of vanishing.
Division before specialization gives
\[
 \left.\frac{\mathcal O_\varepsilon}{\varepsilon^2}
 \right|_{\varepsilon=0}
 =q^{36}r^5d^{12}y^{20}(q^4r^2-1)L\ne0.
\]
This is the same scalar equation, up to a nonzero rational factor,
as the 1998 shift-and-eliminate procedure executed with the
correct shifts, and its published 2001 counterpart.
For sufficiently small \(\varepsilon\ne0\), cyclic rank is
four and both scalar constructions have nonzero fourth-order
leading coefficients. Thus Proposition~\ref{prop:known}'s domain
is nonempty.
\end{proposition}

\section{What the comparison establishes}

The preprint emphasizes generality and explicit coefficients
\cite[pp.~3--4]{MMMRP2026}. Exact translation removes
\(\omega\); the fourth-order Hahn antecedent is
\cite[Theorem~6.5, p.~89, Eq.~(6.59); proof, p.~90]{F1998}.
Different formulas give the same monic operator under
Proposition~\ref{prop:known}.
The published general \(D_q\)-Laguerre--Hahn construction
already uses a \(3\times3\) determinant:
\cite[Theorem~1, p.~456, Eqs.~(45)--(46)]{FHR2001}.
The paper is reference~[10] in the preprint; its own bibliography
lists the thesis as~[15].
The preprint groups [8]--[10] with results for perturbed classical
families \cite[pp.~3--4]{MMMRP2026};
reference~[10] already treats the general class.

For Theorem~2, allow \(B=0\) as in
Proposition~\ref{prop:degrees}. The old first-order relations
\cite[Theorem~6.2, p.~82, Eqs.~(6.29)--(6.32)]{F1998}
give the following explicit verification. Put
\(f=P_{n+1}\), \(g=P_n\), \(\theta=\Theta_n\),
\(\alpha=a+sC_1\), and use \(p,\ell\) from
\eqref{eq:matrices}. At \(n=0\), use \(p=D\),
\(\ell=C_1+(y-\beta_0)D\), and \(R_{-1}=0\).
The first-order residuals are
\begin{equation}
 \begin{aligned}
 E_1&=\alpha\Delta f-(L_n-C_1)f-\theta g+B(TR_n)=0,\\[7pt]
 E_2&=\alpha\Delta g-(\ell-C_1)g+pf+B(TR_{n-1})=0.
 \end{aligned}
 \label{eq:second-order-residuals}
\end{equation}
The left side minus the right side of (56)
in~\cite{MMMRP2026} is exactly
\begin{equation}
 \theta\alpha\Delta E_1+\theta(T\theta)E_2
 -\big\{\alpha\Delta\theta+(T\theta)(\ell-C_1)\big\}E_1.
 \label{eq:second-order-certificate}
\end{equation}
The product rule cancels \(g,\Delta g\);
\(\Delta(TR_n)=qT(\Delta R_n)\) supplies the last coefficient's
factor \(q\). Thus all six coefficients are recovered without
division, including degenerate cases. The old (6.36) and the
published (36), pp.~454--455, of~\cite{FHR2001} retain
different auxiliaries; no literal identity with their tuples is claimed.

Theorem~3 specializes
\cite[Theorem~3.1(i)\(\Rightarrow\)(iii), p.~49,
Eqs.~(3.47), (3.54)--(3.55)]{F1998}
to the original polynomials; its weight and boundary hypotheses give the functional Pearson
identity. For \(N=n+1\), put \(a_2=[y^2]A(c+y)\),
\(\psi_1=[y]\Psi(c+y)\), where \(A,\Psi\) are its classical
coefficient pair. Monicity fixes the eigenvalue by the leading term:
\begin{equation}
 \lambda_N=-[N]_q\big\{\psi_1+a_2q^{1-N}[N-1]_q\big\}
 =\frac{b_{N-1}}q,\quad N\geq1.
 \label{eq:classical-eigenvalue}
\end{equation}
The 1998 expression (3.55) is identical, since
\([N-1]_{1/q}=q^{2-N}[N-1]_q\) and \(A''/(2q)=a_2/q\).
The coefficient pair and
the eigenvalue normalization both coincide.
Both are also printed in \cite[p.~460, Eqs.~(58)--(59)]{FHR2001}.

Section~\ref{sec:bounds} derives Theorem~1's bounds from valid
1998 identities. On the displayed inputs, (33)'s \(k-1\) becomes
the exact degree zero for arbitrarily large \(k\);
\eqref{eq:top-coefficients} gives the leading-term cancellation criteria.
The product identities also have a symmetric \(D,S\) form
on admissible maps~\cite[Definitions~6.3--6.4, p.~32]{C2026}.
Extensions of the Pad\'e bounds and of the scalar comparison
to that setting are outside the scope of this note.

Within the scopes checked above, Theorems~1--3 follow from the
1998 identities and classical result, while Theorem~4 gives the
same monic scalar operator on Proposition~\ref{prop:known}'s domain.
An explicit reformulation does not by itself yield a new scalar
equation on this domain. Deriving degree bounds from old identities
does not establish their prior publication or exclude corrective value.
Expository or computational gains, and additional exceptional coverage,
require their own assessment.

Proposition~\ref{prop:zero-normalization} makes the defect concrete:
on a positive regular family satisfying the stated hypotheses of
Theorem~4, all five prescribed coefficients vanish. The formula
then imposes no condition. Along its explicit regular perturbation,
division by \(\varepsilon^2\) before
specialization recovers a nonzero equation on that family,
proportional to the old scalar equation. A general fourth-order
construction must justify its effective order and specify how its
exceptional specializations are treated.

\section*{Acknowledgements}
The authors acknowledge financial support from the Centre for
Mathematics of the University of Coimbra (CMUC), funded by the
Portuguese Foundation for Science and Technology (FCT), under the
projects UID/00324/2025 and UID/PRR/00324/2025.
The first author acknowledges financial support from the FCT
under grant
2022.00143.\allowbreak{}CEECIND/\allowbreak{}CP1714/\allowbreak{}CT0002.
The second author acknowledges financial support from the FCT
under grant \url{https://doi.org/10.54499/UI.BD.154694.2023}.

The initial draft was produced entirely by \textsf{ChatGPT~6},
in less than fifteen minutes, from photographs of the
\textquotedblleft authors\textquotedblright{}' rough handwritten
notes and calculations containing typos and sign errors.
The AI identified and corrected errors, improved the presentation,
and, on request, reproduced the style of the
\textquotedblleft authors\textquotedblright{}' previous papers.
We subsequently verified the mathematics and assume full
responsibility for every statement here.

\end{document}